%% file: Bourbaki-Viale.tex
\documentclass[brochure,english,12pt]{smfbourbaki}

\usepackage[T1]{fontenc}
\usepackage{lmodern,amssymb,bm}

\usepackage{amsmath, amsthm,amstext,stmaryrd} 
\usepackage{adjustbox}
\usepackage{xcolor}
\usepackage{colortbl}
\usepackage{paracol}
\usepackage{array}   
\newcolumntype{L}{>{$}l<{$}}
\newcolumntype{R}{>{$}r<{$}}
\newcolumntype{C}{>{$}c<{$}}
\newcolumntype{a}{>{\columncolor{LightRed}}l}
\newcolumntype{b}{>{\columncolor{LightCyan}}l}

\definecolor{LightRed}{RGB}{252,160,140}
\definecolor{LightBlue}{RGB}{140,186,252}
\definecolor{Gray}{gray}{0.85}
\definecolor{LightCyan}{rgb}{0.88,1,1}

\usepackage{babel}
\usepackage{graphicx}

\usepackage[utf8]{inputenc}

\usepackage[colorlinks=true, linkcolor=blue, citecolor=red,
urlcolor=blue]{hyperref}

\usepackage[
backend=bibtex,
style=authoryear, 
citestyle=authoryear-comp,
maxnames=7,
sortcites=false 
]{biblatex}
\usepackage{csquotes}

\DefineBibliographyExtras{french}{\restorecommand\mkbibnamefamily}

\DeclareDelimFormat{nameyeardelim}{\addcomma\space}

\DeclareNameAlias{sortname}{given-family}

\renewcommand{\bibnamedash}{\leavevmode\raise3pt\hbox to3em{\hrulefill}\space}

\AtEveryBibitem{%
  \clearfield{issn} 
  \clearfield{isbn} 
  \clearfield{doi} 
  \clearlist{language} 
  \ifentrytype{online}{}{
    \clearfield{url}
  }
}

\renewbibmacro{in:}{%
    \ifentrytype{article}{}{\printtext{\bibstring{in}\intitlepunct}}}

\DeclareFieldFormat[article,periodical,inreference]{number}{\mkbibparens{#1}}
\DeclareFieldFormat[article,periodical,inreference]{volume}{\mkbibbold{#1}}
\renewbibmacro*{volume+number+eid}{%
    \printfield{volume}%
    \setunit*{\addthinspace}
    \printfield{number}%
    \setunit{\addcomma\space}%
    \printfield{eid}}

\DeclareFieldFormat[article,inbook,incollection]{title}{\enquote{#1}\addcomma} 

\date{Avril 2023}
\bbkannee{75\textsuperscript{e} année, 2022--2023}  
\bbknumero{1207}                                      

\title{Strong forcing axioms and the continuum problem}
\subtitle{following Asper\'o's and Schindler's proof that $\mathbf{MM}^{++}$ implies Woodin's Axiom $(*)$}

\author{Matteo Viale}
\address{Dipartimento di Matematica\\ Università di Torino\\ Via Carlo Alberto 10 - 10125\\
Torino - Italie}
\email{matteo.viale@unito.it}

\input{macros}

\begin{document}

\maketitle

\section*{Introduction}
This note addresses the continuum problem, taking advantage of the breakthrough mentioned in the subtitle, and relating it to many recent advances occurring in set theory.\footnote{The author acknowledges support from the project:
\emph{PRIN 2017-2017NWTM8R
Mathematical Logic: models, sets, computability} and from GNSAGA.} 
We try to the best of our possibilities to make our presentation  self-contained and accessible to a general mathematical audience.\footnote{Surveys on the topic complementing this note are (among an ample list) \cite{Bag2005,KOE10,VIAVENRSL,WoodinCH1,WoodinCH2}.}

Let us start by stating Asper\'o's and Schindler's result:

\begin{theo}[\cite{ASPSCH(*)}]
Assume $\mathbf{MM}^{++}$ holds. Then Woodin's axiom $(*)$ holds as well.
\end{theo} 

We will address the following three questions:
\begin{itemize}
\item What is the axiom $\mathbf{MM}^{++}$?
\item What is Woodin's axiom $(*)$?
\item What is the bearing of Asper\'o's and Schindler's result on the continuum problem, and why their result is regarderd as a major breakthrough in the set theoretic community?
\end{itemize}
We give rightaway a spoiler of the type of answers we sketch for the above questions.

We have two major approaches to produce witnesses $x$ of certain mathematical properties $P(x)$.

A topological approach is exemplified for example by Baire's category theorem:  given a compact Hausdorff topological space $X$ one can find a ``generic'' point $x\in X$ satisfying a certain topological property $P(x)$ by showing that $P(x)$ can fail only on a ``small'' (more precisely meager) set of points of $X$.

An algebraic approach is exemplified by the construction of algebraic numbers: one takes a set of Diophantine equations $P(\vec{x})$ which are not jointly inconsistent, and builds abstractly a formal solution in the ring $\mathbb{Q}(\vec{x})/_{P(\vec{x})}$.

Duality theorems connect the algebraic point of view to the geometric one, for example Hilbert's Nullstellensatz relates solutions of irreducible sets of Diophantine equations to generic points of algebraic varieties.

We will outline that Woodin's axiom $(*)$ provides an ``algebraic approach'' to the construction of set theoretic witnesses for ``elementary'' set theoretic properties, $\mathbf{MM}^{++}$ a ``geometric approach'', and Asper\'o's and Schindler's result connects these two perspectives.

We plan to do this while gently introducing the reader to the fundamental concepts of set theory. 

The note is structured as follows:
\begin{itemize}
\item \ref{sec:bassetth} is a brief review of the basic results of set theory with a focus on its historical development and on the topological complexity of sets of reals witnessing the failure of the continuum hypothesis.
\item In \ref{sec:goeprog} we quote some of G\"{o}del's thoughts on the continuum problem and on the ontology of mathematical entities.
\item \ref{sec:larcard} gives a brief overview of (the use in mathematics of) large cardinal axioms.
\item In \ref{sec:forax} we introduce forcing axioms with a  focus on their topological presentations, while giving a precise formulation of the axiom $\mathbf{MM}^{++}$. We also list some of the major undecidable problems which get a solution assuming this axiom, among which the continuum problem.
\item \ref{sec:forc} is a small interlude giving a very fast insight on the forcing method, while relating it to the notions of sheaf and of Grothendieck topos.
\item \ref{sec:modcomp} revolves about the notion of algebraic closure. In particular we outline how Robinson's notion of model companionship gives the means to transfer the concept of ``algebraic closure'' developed for rings to a variety of other mathematical theories.
\item \ref{sec:algclsetth} discusses what is the right language in which set theory should be axiomatized in order to unfold its ``algebraic closure'' properties. 
\item \ref{sec:algclsont} relates Woodin's generic absoluteness results for second order number theory to properties of algebraic closure for the initial fragment of the universe of sets given by $H_{\aleph_1}$.
\item \ref{sec:wooax*}
brings to light why Woodin's axiom $(*)$ can be regarded as an axiom of ``algebraic closure'' for the larger initial fragment of set theory given by $H_{\aleph_2}$.
Putting everything together we conclude by showing why Asper\'o's and Schindler's result establish a natural correspondence between the geometric approach and the algebraic approach to forcing axioms.
\end{itemize}

I thank Alberto Albano, David Asper\'o, Vivina Barutello, Rapha\"{e}l Carroy, Ralf Schindler for many helpful comments on the previous drafts of this manuscript. Many thanks to Nicolas Bourbaki for the invitation and the precious editorial support in the preparation and revision of this work.

\section{Basics of set theory}\label{sec:bassetth}

Set theory deals with the properties of sets (the ``manageable'' mathematical objects) and classes (the ``not so manageable'' entities).\footnote{We refer the reader to \cite{KUNEN,JECHST,MONK} for a systematic treatment of the topic. The reader familiar with set theory can skim through or just skip this section.}

\subsection{Axioms}
The axioms of set theory can be split in three types
(as is the case for many other mathematical theories):
\begin{itemize}
\item \textbf{Universal axioms} which establish properties valid for all sets;
\item \textbf{Existence axioms} which establish the existence of certain sets;
\item \textbf{Construction principles} which allow for the construction of new sets from ones which are already known to exist.
\end{itemize}

We present the axiomatization of set theory by Morse--Kelley $\bool{MK}$ with sets and classes.
Its axioms are distributed in the three categories as follows:
\begin{description}
\item[Universal axioms]
\emph{}
\vspace{0,1cm}
\begin{itemize}
\item \textbf{Extensionality:} 
Two classes (or sets) are equal if they have exactly the same elements.

\item \textbf{Comprehension (a):}  
 Every class (or set) is a subset of $V$, the 
(proper) class whose elements are exactly the sets.

(a \textbf{proper class} is a class which is not a set, a \textbf{set} is a class which belongs to $V$).

\item \textbf{Foundation:} 
There is no infinite sequence $\ap{x_n:n\in\mathbb{N}}$ such that $x_{n+1}\in x_n$ for all $n$.
\end{itemize}

\vspace{0,2cm}
\item[Existence axioms]
\emph{}
\vspace{0,1cm}
\begin{itemize}
\item \textbf{Infinity:} $\emptyset$ and $\mathbb{N}$ are sets.
\end{itemize}

\vspace{0,2cm}
\item[Weak construction principles]
\emph{}
\vspace{0,1cm}

\begin{itemize}
\item \textbf{Union, Pair, Product:} If $X,Y$ are sets, so are $X\cup Y$, $\bp{X,Y}$, $X\times Y$.

\item \textbf{Separation:} If $P$ is a class and $X$ is a set, $P\cap X$ is a set.
\end{itemize}

\vspace{0,2cm}
\item[Strong construction principles]
\emph{}
\vspace{0,1cm}

\begin{itemize}
\item \textbf{Comprehension (b):} For every property $\psi(x)$, 
\(
P_\psi=\bp{a\in V: \psi(a)}
\)
is a class.
\item \textbf{Replacement:} If $F$ is a class function and $X\subseteq\dom(F)$ is a set, the pointwise image $F[X]$ of $X$ under $F$ is a set.

\item \textbf{Powerset:} If $X$ is a set, so is the class
\(
\pow{X}=\bp{Y: Y\subseteq X}.
\)

\item \textbf{Global Choice:} For all classes $C=\bp{X_i:i\in I}$ of non-empty sets $X_i$, 
$\prod_{i\in I}X_i$ (the family of functions $F$ with domain $I$ and such that $F(i)\in X_i$ for all $i\in I$) is non-empty.
\end{itemize}
\end{description}
Some comments: 
\begin{itemize}
\item
By Foundation $V$ cannot be a set else $\ap{x_n:n\in\mathbb{N}}$ with each $x_n$ constantly assigned to $V$ defines a decreasing $\in$-chain.\footnote{\emph{$V$ is not a set} can also be proved without Foundation. Set theorists need foundation in order to infer that the notion of well-foundedness is an elementary set theoretic property (more precisely it is a provably $\Delta_1$-property).} 
\item
Many of the objects of interest in mathematics are proper classes, for example the family of groups, or the family of topological spaces. More generally for a given (first order) theory $T$, the family of stuctures which satisfy the axioms of $T$ is a proper class (and exists in view of Comprehension (b)).
There are delicate ontological issues related to the notion of proper class, but they are foreign to almost all domains of mathematics, with the notable exceptions of category theory and set theory.
\item
It is convenient for natural numbers to distinguish their ordinal type (which confronts them according to which of these numbers ``comes first'') from their cardinal type (which assigns to each natural number $n$ the family of sets which have exactly $n$ elements). When dealing with arbitrary sets, their ordinal type may not be defined, while the cardinal type always is. 
Von Neumann devised a simple trick to represent the finite ordinal types. One can inductively define the natural number $n$ as the set $\bp{0,\dots,n-1}$ (i.e.\ $0=\emptyset$, $1=\bp{\emptyset}$, $2=\bp{\emptyset,\bp{\emptyset}}$, 
\dots).\footnote{The transfinite ordinal types (or the Von Neumann ordinals) are those (possibly infinite) sets $\alpha$ which are linearly ordered by $\in$ and are transitive (i.e.\ such that when $x\in y\in\alpha$, we have that $x\in\alpha$ as well). The proper class of Von Neumann ordinals is linearly well-ordered by $\in$. One can check that the natural numbers are the finite Von Neumann ordinals and that $\mathbb{N}$ is the first infinite Von Neumann ordinal.}
\item
Set theoretic construction principles are of two sorts: the simple (or weak) ones are for example those bringing from sets $X,Y$
to sets $X\cup Y$, $\bp{X,Y}$, $X\times Y$, or from set $X$ and class $P$ to the set $P\cap X$; the strong ones are the power-set axiom, the replacement axiom, and the axiom of choice. 
Let us discuss briefly the role of such axioms in the development of routine mathematics.
\begin{description}
\item[Weak construction principles]
The integers and rationals can be constructed from the naturals using only weak construction principles: 
$\mathbb{Z}$ can be seen as the subset of $\mathbb{N}\times\bp{0,1}$ which assigns the positive integers to the ordered pairs with second coordinate $0$ and the negative ones to those pairs with second coordinate $1$ (paying attention to the double counting of $0$ as $(0,0)$ and $(0,1)$);  
$\mathbb{Q}$ can be seen as the subset of $\mathbb{Z}\times(\mathbb{N}\setminus\bp{0})$ given by ordered pairs which are coprime.

\item[Powerset axiom]
In order to build the reals from the rationals, one needs this axiom: $\mathbb{R}$ is the subset of $\pow{\mathbb{Q}}$ given by Dedekind cuts.

\item[Replacement axiom]

%
%

An adequate development of set theory requires it: consider the function $F$ on $\mathbb{N}$ given by $F(0)=\mathbb{N}$, $F(n+1)=\pow{F(n)}$.
Without replacement it cannot be proved that $F$ (or even the image of $F$) is a set, it might only be a proper class.


\item[Choice]
Choice also has a special status in ordinary mathematics, and many mathematicians feel uneasy about it. However Choice is unavoidable: it is essential in the proofs of the Hahn--Banach theorem, of the existence of a base for infinite-dimensional vector spaces, or of the existence of a maximal ideal on a ring,\dots Even the equivalence of sequential continuity and topological continuity for real valued functions requires it: if $f:\mathbb{R}\to\mathbb{R}$ is not continuous at $x$, there is $\epsilon>0$ such that for each $n$ one can find $x_n$ so that $|x_n-x|<1/n$ and $|f(x_n)-f(x)|>\epsilon$.
The sequence $(x_n)_n$ is (and in most cases can only be) defined appealing to (countable) Choice.
\end{description}
\end{itemize}

\subsection{Cardinal arithmetic}

Let us now develop arithmetic on cardinal types.

Given sets $X$, $Y$
\begin{itemize}
\item The cardinal $|X|$ is the (proper) class $\bp{Y: \, \exists \,f:X\to Y \text{ bijection}}$;

\item $|X|\leq|Y|$ if and only if \emph{there is $f:X\to Y$ injection} if and only if \emph{there is $g:Y\to X$ surjection};\footnote{The last equivalence gives an alternative definition of the axiom of choice.}

\item $|X|<|Y|$ iff $|X|\leq|Y|$ and $|X|\neq|Y|$;

\item $|X|+|Y|=|(X\times\bp{0})\cup (Y\times\bp{1})|$ (the size of the disjoint union of $X$ and $Y$);
\item $|X|\cdot |Y|=|X\times Y|$ (the size of the product of $X$ and $Y$);
\item $|X|^{|Y|}=|X^Y|$ where 
\[
X^Y=\bp{f: \text{ $f$ is a function with domain $Y$ and range $X$}}.
\]
\end{itemize}
Note that each cardinality class on a set $X$ is a proper class if it contains a non-empty element.
The equivalence class of $\emptyset$ is $\bp{\emptyset}$.
It can also be shown that $C$ is a proper class if and only if there is a surjection of $C$ onto $V$, and 
$X$ is a set if and only if there is no surjection of $X$ onto $V$.

The finite sets are those in the equivalence class of some $n\in\mathbb{N}$. For finite sets $X,Y$ of size respectively $m,n$:
\begin{itemize}
\item $|X|=|Y|$ if and only if $n=m$ if and only if $|X|\leq|Y|$ and $|Y|\leq|X|$;
\item $|X|+|Y|$ is (in  the equivalence class of) $m+n$;
\item $|X|\cdot|Y|$ is (in  the equivalence class of) $m\cdot n$;
\item $|X|^{|Y|}$ is (in  the equivalence class of) $m^n$;
\end{itemize}
i.e.\ the map $n\mapsto |n|$ defines an embedding of the stucture 
$(\mathbb{N},<,\cdot,+,n\mapsto 2^n,0,1)$ in the family of cardinals endowed by the above operations.

Basic questions were then asked for the infinite cardinals:
\begin{enumerate}
\item Assume $|X|\leq|Y|$ and $|Y|\leq|X|$ holds for infinite sets $X,Y$, do we have also in this case that $|X|=|Y|$?
\item What is the structure of the order on infinite cardinals given by $<$?
\item Do the arithmetic of cardinals obey the associative, commutative, distributive laws which hold for the naturals?
\item What is the computation table for the arithmetic operations of sum and product on cardinals, when one of the factors is an infinite cardinal?
\item What about the computation table for the exponential map $|X|\mapsto 2^{|X|}$ when $X$ is infinite?
\end{enumerate}

Soon most of these questions had a clear cut answer:
\begin{enumerate}
\item For arbitrary sets $X,Y$ we have that $|X|\leq|Y|$ and $|Y|\leq|X|$ if and only if $|X|=|Y|$ (Cantor 1887, Bernstein 1897, Dedekind 1898).\footnote{This is non trivial: there is a topological embedding of $[0;1]$ into $(0;1)$ and conversely; there is no continuous bijection of $[0;1]$ onto $(0;1)$; however one can find a Borel bijection.}
\item \label{cardZER}
$\leq$ is a well-order on cardinals (Zermelo+\dots $\sim$ 1904), i.e.\ it is a linear order on cardinals such that for every class $C\neq\emptyset$ there is $\min\bp{|X|: X\in C}$.\footnote{
Furthermore in each class $|X|$ there is a least Von Neumann ordinal $\kappa$ which is called the cardinal of $X$ and is the canonical representative of $|X|$.}

\item
The arithmetic operations on cardinals endow this family of a structure of ordered commutative semiring with exponentiation; the natural numbers form an initial segment of this semiring (Cantor - published in 1895, but most likely known earlier).
\item $|X|+|Y|=|X|\cdot|Y|=\max\bp{|X|,|Y|}$ if at least one among $X,Y$ is infinite and both are non-empty (Special cases by Cantor, before 1895, general proofs by Harward 1905, Jourdain 1908, Hausdorff 1914\dots).
\item \label{cardCAN}
For all sets $X$, $|X|<|\pow{X}|=|2|^{|X|}$ (Cantor 1891).\footnote{$x\mapsto\bp{x}$ witnesses $|X|\leq|\pow{X}|$; if $g:\pow{X}\to X$, $Y_g=\bp{y\in Y:y\not\in g(y)}$ witnesses that $g$ is not a surjection.

Regarding the map $(|X|,|Y|)\mapsto |X|^{|Y|}$, the computation of its table can be done in terms of that of the exponential map \parencite{BUK65}. }
\end{enumerate}
We denote cardinals by greek letters $\kappa,\lambda,\theta$. 
\begin{itemize}
\item
$\aleph_0$ is the cardinality of $\mathbb{N}$ the least infinite cardinal;
\item
for $\kappa$ a cardinal $\kappa^+$
 is the cardinality of the least cardinal above $\kappa$ (which exists by \ref{cardZER} and \ref{cardCAN});
 \item
 $\aleph_1=\aleph_0^+$, $\aleph_2=\aleph_1^+$.
 \end{itemize}

\subsection{The continuum problem}
The so called continuum problem (Cantor 1878) remained unsettled:
\begin{quote}
What is the value of $2^{|\mathbb{N}|}$?
\end{quote}

The continuum hypothesis can be equivalently phrased:
\begin{itemize}
\item $2^{\aleph_0}=\aleph_1$;
\item For all $X\subseteq \mathbb{R}$ either $X$ is countable or $X$ has size continuum.
\end{itemize}
The equivalence of the two follows once one observes that $|\mathbb{R}|=|\pow{\mathbb{N}}|$.

\subsection{Progresses on the definable version of the continuum problem}
For decades there were scant progresses on the solution of this problem.
On the positive side one can note proofs that ``simply definable'' subsets of $\mathbb{R}$ cannot witness the negation of the continuum hypothesis (the  so called ``continuum hypothesis for definable sets'').

The \emph{analytic} (or $\Sigma^1_1$) subsets of $\mathbb{R}^n$ are the projections on the first $n$-coordinates of a Borel subset of $\mathbb{R}^{n+k}$ for some natural number $k$. The \emph{coanalytic} (or $\Pi^1_1$) subsets of $\mathbb{R}^n$ are the complements of the analytic subsets.

\begin{itemize}
\item  No \emph{closed} subset of $\mathbb{R}$ is a counterexample to $\bool{CH}$ 
(Cantor 1883).

\item  No \emph{Borel} subset of $\mathbb{R}$ is a counterexample to $\bool{CH}$ 
(Alexandroff 1916, Hausdorff 1917).

\item  No \emph{analytic} subset of $\mathbb{R}$ is a counterexample to $\bool{CH}$ 
(Suslin+Alexandroff 1917).

\item  \emph{Coanalytic} subsets of $\mathbb{R}$ have either size $\aleph_0,\aleph_1,2^{\aleph_0}$ (Luzin--Sierpinski 1917$\approx$).
\end{itemize}

This was as far as the first rounds of mathematicians addressing Cantor's problem could get.

With the introduction of large cardinals further major progresses were obtained on the definable version of $\bool{CH}$.

The \emph{projective} subsets of $\mathbb{R}^n$ are those subsets of $\mathbb{R}^n$ which are $\Sigma^1_m$
(or $\Pi^1_m$) for some~$m$, where $X\subseteq \mathbb{R}^n$ is $\Sigma^1_{m+1}$ if it is the projection of a $\Pi^1_{m}$-subset of  $\mathbb{R}^{n+k}$ and is  a $\Pi^1_{m+1}$-subset if it is the complement of a $\Sigma^1_{m+1}$-set.

Projective sets define a natural family of subsets of $\mathbb{R}^n$ which ought to be topologically simple.

Another natural family has been isolated by \textcite{FENMAGWOO}:

Recall that for a given topological space $(Y,\tau)$, a set $X$ is \emph{nowhere dense} in $Y$ if its complement contains an open dense subset of $Y$, it is \emph{meager} if it is a countable union of nowhere dense subsets, it has the \emph{Baire property} if it has meager symmetric difference with an open set.
\begin{defi}[\cite{FENMAGWOO}]\label{def:UB}
$X\subseteq \mathbb{R}^k$ is \emph{universally Baire} if for all continuous maps $f:Y\to\mathbb{R}^k$
with $Y$ compact Hausdorff, $f^{-1}[X]$ has the Baire property in $Y$.
\end{defi}
Analytic sets are universally Baire, and this family forms a $\sigma$-algebra.

To appreciate the strength of this property, consider $2^{\mathbb{N}}$ when $2$ is given the discrete topology and $2^{\mathbb{N}}$ has the product topology.

The map $\theta:f\mapsto\sum_{i=0}^\infty f(i)/3^{i+1}$ for $f:\mathbb{N}\to\bp{0,2}$ defines a topological embedding of $2^\mathbb{N}$ into $[0;1]$ whose image 
is meager and has Lebesgue measure $0$.
%

Now take a subset $P$ of $2^{\mathbb{N}}$ which does not have the Baire property in $2^{\mathbb{N}}$.

Seen as a subset of $[0;1]$, $\theta[P]$ is meager, hence it has the Baire property in $[0;1]$.
On the other hand $P=\theta^{-1}[\theta[P]]$ does not have the Baire property in $2^{\mathbb{N}}$.
Hence $\theta[P]$ is not universally Baire even if it has the Baire property.

A similar argument based on measure, shows that the non-measurable sets produced by Vitali are not universally Baire. More generally it can be shown that assuming large cardinals a subset of the reals is ``non-pathological'' if and only if it is universally Baire.
The theorem below makes this a sound mathematical assertion.
\begin{theo}
Assume there is a proper class of Woodin cardinals. Then:
\begin{itemize}
\item  No universally Baire subset of $\mathbb{R}$ is a counterexample to $\bool{CH}$ 
\parencite{FENMAGWOO,MARSTE89,DAVIS64}.
\item  Borel sets, analytic sets, projective sets,\dots are all universally Baire 
\parencite{FENMAGWOO}.
\item The universally Baire subsets of $\mathbb{R}$ are determined \parencite{MARSTE89}.
\item The determined subsets of the reals are Lebesgue measurable, have the perfect set property,\dots
(various authors one for each relevant regularity property).
\end{itemize} 
\end{theo}

\subsection{Independence of $\bool{CH}$}
$\bool{CH}$ is independent of the axioms of set theory:
\begin{itemize}
\item
There is a model of the axioms of $\bool{MK}$ where $\bool{CH}$ holds \parencite{GOD47}. 
\item  
There is a model of the axioms of $\bool{MK}$ where $\bool{CH}$ fails \parencite{COH63}.
\item
In the model of the axioms of $\bool{MK}$ where $\bool{CH}$ fails produced by Cohen, this failure can be witnessed by a $\Sigma^1_2$-set.
\end{itemize}


\section{{G}\"{o}del's program} \label{sec:goeprog}
\textcite{GOD47} wrote an influential survey on the continuum problem which has been a source of inspiration for the work of many logicians to come. In my opinion the best summary of its content can be given by quoting a few excerpts from the paper:

\begin{description}
\item[On the independence of the continuum problem (p.~520)]
\emph{Only someone who (like the intuitionist) denies that the concepts and axioms of classical set theory have any meaning (or any well-defined meaning) could be satisfied with such a solution, not someone who believes them to describe some well-determined reality. For in this reality Cantor's conjecture must be either true or false, and its undecidability from the axioms as known today can only mean that these axioms do not contain a complete description of this reality;}
\item[On Large Cardinals (p.~520)]
\emph{For first of all the axioms of set theory by no means form a system closed in itself, but, quite on the contrary, the very concept of set on which they are based suggests their extension by new axioms which assert the existence of still further iterations of the operation ``set of''. These axioms can also be formulated as propositions asserting the existence of very great cardinal numbers or (which is the same) of sets having these cardinal numbers. The simplest of these strong ``axioms of infinity'' assert the existence of inaccessible numbers (and of numbers inaccessible in the stronger sense) $>\aleph_0$.}
 
\item[On success as a criterion to detect new axioms (p.~521)]
\emph{There might exist axioms so abundant in their verifiable consequences, shedding so much light upon a whole discipline, and furnishing such powerful methods for solving given problems (and even solving them, as far as that is possible, in a constructivistic way) that quite irrespective of their intrinsic necessity they would have to be assumed at least in the same sense as any well established physical theory.} 
\end{description}




\section{Large cardinals} \label{sec:larcard}
Consider the universe of sets construed from the emptyset using all other axioms with the exception of \emph{Infinity}; it can be shown that one ends up with a ``baby'' universe $H_{\aleph_0}$ characterized by the following property: if $X\in H_{\aleph_0}$, $X$ is finite, all the elements~$Y$ of~$X$ are finite, all the elements~$Z$ of some element~$Y$ of~$X$ are finite,\dots
More precisely
\begin{defi}\label{def:Haleph0}
Given a set $X$, 
\begin{itemize}
\item
$\cup^0 X=X$, $\cup^{n+1}X=\cup(\cup^n X)$, 
\item
$\trcl(X)=\bigcup_{n\in\mathbb{N}}(\cup^n X)$ is the \emph{transitive closure} of $X$.
\end{itemize}
$X$ is \emph{hereditarily finite} if $\trcl(X)$ is finite.

$H_{\aleph_0}$ is the collection of hereditarily finite sets.

%

\end{defi}

%
%
%
%
%
One can check that $\pow{H_{\aleph_0}}$ (or -more precisely- the structure $(\pow{H_{\aleph_0}},H_{\aleph_0},\in)$) is a model of all axioms of Morse--Kelley set theory with the exception of the infinity axiom asserting that $\mathbb{N}$ is a set. Indeed in this model $H_{\aleph_0}$ is the class of all sets and $\mathbb{N}\subseteq H_{\aleph_0}$ is a ``proper class''.

In particular the axiom of Infinity is a way to assert that not all sets can be described using the construction principles and starting from the emptyset. $\mathbb{N}$ is an example of such a set.

Large cardinal axioms posit the existence of sets which cannot be described from $\mathbb{N}$ using the
construction principles encoded in the axioms of $\bool{MK}$. 
We already mentioned from {G}\"{o}del the axiom stating the existence of inaccessible cardinals, i.e.\ cardinals behaving like $|\mathbb{N}|$, but of larger size: $\kappa$ is inaccessible if and only if $\pow{H_\kappa}$ is a model of Morse--Kelley set theory whose universe of sets is $H_\kappa$,\footnote{Equivalently $\kappa$ is inaccessible if it is regular and $(H_\kappa,\in)\models\ZFC$.} where:
\begin{defi}\label{def:Hkappa}
Given a cardinal $\lambda$, a set $X$ is \emph{hereditarily of size less than $\lambda$} if $\trcl(X)$ has size less than $\lambda$.

$H_\lambda$ is the set of all sets which are hereditarily of size less than $\lambda$.
\end{defi}

The original proof by Wiles of Fermat's last theorem uses the notion of Grothendieck universe; when correctly formalized in set theory, the existence of a Grothendieck universe is equivalent to the existence of an inaccessible cardinal. Grothendieck's theory of universes finds its natural formulation in (it is actually equivalent to) set theory enriched with the axiom stating the existence of a proper class of inaccessible cardinals \parencite{MCLAR10}.

Wiles' proof provides evidence grounded on {G}\"{o}del's criterion of success  for the adoption of large cardinal axioms. There are plenty of (less celebrated) such cases.
For example:

\begin{defi}[Vopenka's principle $\bool{VP}$]
For every \emph{proper class} of \textbf{directed graphs with no loops}, there are two members of the class with a homomorphism between them.
\end{defi}
This is an equivalent formulation of $\bool{VP}$ by \textcite{ADAROS94}. $\bool{VP}$ is a tool which category theorists employ successfully, see for example \cite{BAGCASMATROS15,CASSCESMI05,ROSTHO03}.

$\bool{VP}$ entails the existence of all large cardinal axioms one may require in any of the results presented elsewhere in this note; for example if $\bool{VP}$ holds there is a proper class of supercompact cardinals; furthermore any supercompact cardinal is also Woodin.
Henceforth assuming $\bool{VP}$ no counterexample to $\bool{CH}$ can be universally Baire, the projective sets of reals are universally Baire, the universally Baire sets are determined.

Later on we will delve more on the effects of large cardinals on second order numebr theory and on the family of topologically simple sets of reals.

\section{Forcing axioms}\label{sec:forax}

Loosely speaking forcing axioms try to encapsulate the idea that the powerset of some set $X$ is ``as thick as possible''.
Forcing axioms for $X$ can be divided in two categories:

\begin{description}

\item[Axioms of topological maximality] They can be reformulated as strong forms of Baire's category theorem, are inspired by the notion of generic point, include
$\MM^{++}$ among their instantiations.

\item[Axioms of algebraic maximality] They assert the closure of  $\pow{X}$ under a variety of set theoretic operations, are inspired by the notion of 
algebraic closure,  include Woodin's axiom $(*)$ among their instationations. 

\end{description}

We remark the following:
\begin{itemize}
\item
$\MM^{++}$ and $(*)$ are forcing axioms for $X=\aleph_1$.
\item
Baire's category theorem is a ``topological'' forcing axiom for $X=\mathbb{N}=\aleph_0$.

\item
Large cardinals entail ``algebraic'' forcing axioms for $X=\mathbb{N}=\aleph_0$.

\end{itemize}

\subsection{Topological maximality and Martin's maximum}

Recall Baire's category theorem:
\begin{quote}
Let $(X,\tau)$ be a compact Hausdorff space and $\bp{D_n:n\in\mathbb{N}}$ be a countable family of dense open subsets of $X$. Then $\bigcap_n D_n$ is dense in $X$.
\end{quote}

Let us parametrize the conclusion in all cardinals $\kappa$ rather than just $\aleph_0$:
\begin{defi}
Let $\kappa$ be an infinite cardinal and $(X,\tau)$ a topological space. 

$\bool{FA}_{\kappa}(X,\tau)$ holds if $\bigcap_{i\in\kappa}D_i$ is \emph{dense} in $X$ for all $\bp{D_i:\, i \in\kappa}$ family of dense open subsets of $X$.
\end{defi}

For $\kappa>\aleph_0$ not all $(X,\tau)$ compact Hausdorff satisfy $\FA_\kappa(X,\tau)$. For example:
\begin{quote}
Let $Y$ be an \emph{uncountable set} and $(X,\tau)$ be the Stone-\v{C}ech compactification of the product space $Y^{\mathbb{N}}$ where $Y$ is endowed of the discrete topology. Then 
$\FA_{\aleph_1}(X,\tau)$ fails.
\end{quote}
It can also be shown that $\FA_{\aleph_1}(X,\tau)$ holds for certain compact Hausdorff spaces $(X,\tau)$.\footnote{Actually an equivalent formulation of the axiom of choice states that $\FA_\kappa(X,\tau)$ holds for any compact Hausdorff space $(X,\tau)$ such that $\tau$ admits a base $P$ of non-empty sets with the property that $(P,\subseteq)$ is a $<\kappa$-closed forcing (see for example \cite{VIAUSAX}).}

Abraham isolated a necessary condition on $(X,\tau)$ so that 
$\bool{FA}_{\aleph_1}(X,\tau)$ is not inconsistent.

\begin{prop}[Abraham]
Assume $(X,\tau)$ is a compact Hausdorff space which is not\footnote{See Def. \ref{def:SSP}; however now we do not need to know what $\SSP$ means for now.} $\bool{SSP}$.
Then $\FA_{\aleph_1}(X,\tau)$ fails.
\end{prop}

\textcite{FORMAGSHE} showed that it can also be a sufficient condition:

\begin{defi}[\cite{FORMAGSHE}]
Martin's maximum $\MM$ $\equiv$\emph{ $\FA_{\aleph_1}(X,\tau)$ holds for all compact Hausdorff spaces $(X,\tau)$ which are $\bool{SSP}$.}
\end{defi}

\begin{theo}[\cite{FORMAGSHE}]
Assume there exists a proper class of supercompact cardinals. Then there is a model of $\bool{MK}$ where
$\MM$ holds and there is a proper class of supercompact cardinals.
\end{theo}

In particular the theorem establishes that there is a model of set theory enriched with large cardinal axioms such that
\begin{quote}
\emph{
$\FA_{\aleph_1}(X,\tau)$ holds for all compact Hausdorff spaces $(X,\tau)$  for which it is not impossible.}
\end{quote}

\subsection{What are the mathematical consequences of forcing axioms?}
We collect here some of the major applications of forcing axioms in set theory and in  other domains
of mathematics.

Assume Martin's maximum holds. Then:
\begin{itemize}
\item
$\bool{CH}$ is false and
the continuum is the second uncountable cardinal, i.e.\
 $2^{\aleph_0}=\aleph_2$ \parencite{FORMAGSHE}.

\item
Whitehead's conjecture on free groups is false,
(i.e.\ there are uncountable Whitehead groups which are not free) \parencite{SHE74}.

\item Kaplansky's conjecture on Banach algebras holds (i.e. every algebra homomorphism from the Banach algebra $C(X)$ -where $X$ is compact Hausdorff- into any other Banach algebra is necessarily continuous) 
(Woodin and Solovay, unpublished).\footnote{A complete proof from published sources can be obtained combining results in \cite[Chapter 3]{DALWOO87} with \cite[Thm. 7.7]{TOD89}.}

\item
There are five uncountable linear orders such that any uncountable linear order contains an isomorphic copy of one of them \parencite{moore.basis}

\item
All automorphisms of the Calkin algebra are inner \parencite{FAR11}. 

\item \dots
\end{itemize}

All these conclusions are independent of $\bool{MK}+$\emph{Vopenka's principle} (or any other large cardinal assumption).

\subsection{$\bool{MM}^{++}$}
$\bool{MM}^{++}$ is a natural technical strengthening of $\bool{MM}$.
For the sake of completeness we state one of its possible definitions; we caution the reader that it will make sense only for those familiar with the forcing method and with the basic theory of $\pow{\aleph_1}$.
Those not willing to delve into set theoretic technicalities may skip to the next section. This axiom also appears in \textcite{FORMAGSHE}, as well as all the results of this section, unless otherwise specified.

\begin{defi}
Given a complete boolean algebra $\bool{B}$, $\dot{S}\in V^{\bool{B}}$ is a $\bool{B}$-name for a stationary subset of $\aleph_1$, if 
\(
\Qp{\dot{S}\text{ is a stationary subset of }\dot{\aleph}_1}_{\bool{B}}=1_{\bool{B}}.
\)

Given $G$ ultrafilter on $\bool{B}$, we let 
\(
\dot{S}_G=\bp{\alpha<\aleph_1: \Qp{\check{\alpha}\in\dot{G}}_{\bool{B}}\in G}.
\)

Given $C$ club subset of $\aleph_1$, we let 
\(
D_{C,\dot{S}}=\bp{G\in\St(\bool{B}): \dot{S}_G\cap C\neq\emptyset}.
\)
\end{defi}

Note that  $D_{C,\dot{S}}$ is dense open in $\St(\bool{B})$ for any $C$ club subset of $\aleph_1$.

Recall that $(X,\tau)$ is a compact extremally disconnected Hausdorff space if and only if it is the Stone space of its algebra of regular open sets. 

\begin{prop}
For a compact Hausdorff space $(X,\tau)$,
$\FA_{\aleph_1}(X,\tau)$ holds if and only if so does $\FA_{\aleph_1}(\St(\bool{B}),\tau_{\bool{B}})$;
where $\bool{B}$ is the algebra of regular open subsets of $X$ and $\tau_{\bool{B}}$ is the compact Hausdorff topology on the Stone space $\St(\bool{B})$ of ultrafilters on $\bool{B}$ induced by the Stone duality applied to $\bool{B}$.
\end{prop}

\begin{defi}
Given a compact extremally disconnected Hausdorff space $(X,\tau)$, we let $\bool{B}$ be the complete 
boolean algebra given by its regular open sets. 

$\FA^{++}_{\aleph_1}(X,\tau)$ holds if for any family $\bp{D_\alpha:\alpha<\aleph_1}$ of dense open subsets of $X$, and any family $\bp{\dot{S}_\alpha:\alpha<\aleph_1}$ of $\bool{B}$-names for stationary subsets of $\aleph_1$, there is a  ultrafilter $G\in X=\St(\bool{B})$ such that
\[
G\in D_\alpha \text{ for all }\alpha<\aleph_1,
\]
and
\[
G\in D_{C,\dot{S}_\alpha} \text{ for all }\alpha<\aleph_1 \text{ and $C$ club subset of $\aleph_1$}.
\]
\end{defi}
Note that $2^{\aleph_1}>\aleph_1$ and there are $2^{\aleph_1}$-many clubs $C$ on the cardinal $\aleph_1$, each one producing a different set $D_{C,\dot{S}_\alpha}$ for each $\alpha<\aleph_1$. Therefore $\FA^{++}_{\aleph_1}(X,\tau)$ requires the existence of points on $X$ meeting families of dense open sets of size much larger than~$\aleph_1$.

\begin{defi}
$\bool{MM}^{++}$ holds if $\FA^{++}_{\aleph_1}(X,\tau)$ holds for all compact extremally disconnected Hausdorff spaces which are $\bool{SSP}$.
\end{defi}

It is known (for example by a combination of \cite{ASPSCH(*),larson-dwo}) that $\bool{MM}^{++}$ is an axiom strictly stronger than $\bool{MM}$. However the standard proof of the consistency of Martin's maximum produces a model of $\bool{MM}^{++}$, and we have the following:

\begin{theo}
Assume there is a proper class of supercompact cardinals. Then there is a model of $\bool{MK}$ with a proper class of supercompact cardinals where
$\MM^{++}$ holds.
\end{theo}

\section{Forcing}\label{sec:forc}
We give in this section some basic information on the forcing method and its relation with the notion of Grothendieck Topos. The content of this section is not needed in the subsequent parts of this note.

We recall that Stone duality identifies complete boolean algebras with compact Hausdorff extremally disconnected spaces (a space is extremally disconnected if the closure of an open set is open).

Below we organized a text sketching on its left-side column the forcing procedure according to a set theorist, and highlighting on the right-side column the corresponding steps viewed with the lenses of a category theorist.
The reader can skim through the left-side text, then through the right side text, and finally compare the paragraphs of the two texts with same alignment. 

\vspace{1cm}

\begin{paracol}{2}
\centerline{\textbf{SET THEORIST}}
\switchcolumn
\centerline{\textbf{CATEGORY THEORIST}}
\switchcolumn

\smallskip

\switchcolumn

\smallskip

\switchcolumn
\noindent
Given the complete boolean algebra $\bool{B}$,

\switchcolumn
\noindent
Given the compact extremally disconnected Hausdorff space $\St(\bool{B})$,
\end{paracol}

\medskip

\centerline{one forms}

\medskip

\begin{paracol}{2}
\noindent the boolean valued model $V^{\bool{B}}$ \\
by Cohen/Scott--Solovay--Vopenka;
\switchcolumn
\noindent the topos $\bool{Sh}(\St(\bool{B}),\bool{CompHaus})$, given by topological sheaves on $\St(\bool{B})$ with target a compact Hausdorff space;
\end{paracol}

\medskip

\centerline{one chooses}

\medskip

\begin{paracol}{2}
\noindent a $V$-generic ultrafilter $G$ on $\bool{B}$,
\switchcolumn
\noindent a generic point of $\St(\bool{B})$ belonging to all dense open subset of $\St(\bool{B})$,
\end{paracol}

\medskip

\centerline{and one obtains}

\medskip

\begin{paracol}{2}
\noindent the $\bool{MK}$-model $V[G]$ generic extension of $V$ by $G$.
\switchcolumn
\noindent the topos $\bool{Sh}(\St(\bool{B}),\bool{CompHaus})/_G$ of stalks at $G$ of the sheaves in\\
$\bool{Sh}(\St(\bool{B}),\bool{CompHaus})$.
\end{paracol}

\medskip

\medskip

\medskip

\centerline{The properties of}

\medskip

\begin{paracol}{2}
\noindent 
$V[G]$ depend mainly on $\bool{B}$
\switchcolumn
\noindent
$\bool{Sh}(\St(\bool{B}),\bool{CompHaus})/_G$ depend mainly on $\St(\bool{B})$
\end{paracol}

\medskip

\centerline{and minimally on $G$.\footnote{Note that generic points do not exist for atomless boolean algebras, however  it can be of help for our intuition to work under the assumption that such points exist and that -on average- the points of $\St(\bool{B})$ behave like generic points. A strictly formal approach to the topic can eliminate the use of generic points in the forcing analysis of the properties of $V^\bool{B}$ (or equivalently of $\bool{Sh}(\St(\bool{B}),\bool{CompHaus})$).}}

We give some glimpses on how forcing assigns truth values to logical properties in a structure of the form
$\mathcal{F}(X)/_G$ when $\mathcal{F}$ is a sheaf on a compact extremally disconnected Hausdorff space $X$, and $G$ is a point of this space.
The machinery of forcing generalizes to the topos of such sheaves this procedure.

Consider the space $L^{\infty+}(\mathbb{R})$ of measurable functions $f:\mathbb{R}\to\mathbb{R}\cup\bp{\infty}$ which take value $\infty$ on a null set. For example $1/x$ is such a function and all essentially bounded measurable functions are also in $L^{\infty+}(\mathbb{R})$. Note that with pointwise multiplication and sum $L^{\infty+}(\mathbb{R})/_\approx$ is a commutative ring (where $\approx$ identifies functions overlapping an a conull set, and $f+g(x)=f\cdot g(x)=0$ if either one of $f(x),g(x)$ is $\infty$).

Let $\bool{MALG}$ be the complete boolean algebra given by Lebesgue measurable sets modulo null sets.
Given a universally Baire (hence Lebesgue measurable) relation $R\subseteq \mathbb{R}^n$ we can lift it to a boolean relation 
\[
R^{\bool{MALG}}:L^{\infty+}(\mathbb{R})^n\to\bool{MALG}
\]
by setting 
\[
R^{\bool{MALG}}(f_1,\dots,f_n)=[\bp{x\in\mathbb{R}: R(f_1(x),\dots,f_n(x))}],
\]
where $[A]$ is the equivalence class in $\bool{MALG}$ of the measurable set $A\subseteq\mathbb{R}$.
For example
\[
\sin(x)<^{\bool{MALG}}\cos(x)=[\bigcup_{n\in\mathbb{Z}}((2n-1)\cdot\pi+\frac{\pi}{4};2n\cdot\pi+\frac{\pi}{4})]
\]
gets a $\bool{MALG}$-value which is neither true nor false.

A key observation is that whenever $R,S\subseteq \mathbb{R}^n$ are universally Baire (Lebesgue measurable), so are
$R\cap S,R\cup S, \mathbb{R}^n\setminus S$. Furthermore (assuming the existence of a proper class of Woodin cardinals)
so is also $\pi_{j}[R]$ where $\pi_j$ is the projection on the $j$-th coordinate of $R$.

This gives that
\begin{align*}
(R\cap S)^{\bool{MALG}}(f_1,\dots,f_n)&\\ =&[\bp{x\in\mathbb{R}: R(f_1(x),\dots,f_n(x))}\cap \bp{x\in\mathbb{R}: S(f_1(x),\dots,f_n(x))}]\\
=&R^{\bool{MALG}}(f_1,\dots,f_n)\wedge S^{\bool{MALG}}(f_1,\dots,f_n)
\end{align*}
Similarly it can be shown (appealing again to universal Baireness) that:
\[
(\mathbb{R}^n\setminus S)^{\bool{MALG}}(f_1,\dots,f_n)=\neg_{\bool{MALG}} S^{\bool{MALG}}(f_1,\dots,f_n),
\]
\[
(R\cup S)^{\bool{MALG}}(f_1,\dots,f_n)=R^{\bool{MALG}}(f_1,\dots,f_n)\vee S^{\bool{MALG}}(f_1,\dots,f_n),
\]
\[
(\pi_n[R])^{\bool{MALG}}(f_1,\dots,f_{n-1})=\bigvee_{\bool{MALG}}
\bp{R^{\bool{MALG}}(f_1,\dots,f_{n-1},g):g\in L^{\infty+}(\mathbb{R})}.
\]
Now when $G$ is a ultrafilter on $\bool{MALG}$, we can pass to the quotient structure
$L^{\infty+}(\mathbb{R})/_G$ given by the classes 
\[
[f]_G=\bp{h: [f=h]\in G}
\]
and relations
\[
R/_G([f_1]_G,\dots,[f_n]_G)
\]
holding when $R^{\bool{MALG}}(f_1,\dots,f_n)\in G$.

One can check (among many things) that the axioms of ordered field holds in $L^{\infty+}(\mathbb{R})/_G$, as $G$ makes a coherent selection of truth values. For example  given $f,g,h\in L^{\infty+}(\mathbb{R})$
it cannot be the case that the three boolean values $(f<^{\bool{MALG}}g), (g<^{\bool{MALG}}h), (h<^{\bool{MALG}}f)$ are all simultaneously in $G$. Similar arguments yield that the relation $<^{\bool{MALG}}/_G$ defines a dense linear order without endpoints on $L^{\infty+}(\mathbb{R})/_G$.

Letting $f\approx g$ iff $[\bp{x:f(x)=g(x)}]=1_{\bool{MALG}}$, $L^{\infty+}(\mathbb{R})/_\approx$ is the space of global sections with regard to the sheafification according to the dense Grothendieck topology of the presheaf $\mathcal{F}$ which to any
$[A]\in\bool{MALG}$ assigns the family of essentially bounded measurable functions defined on $A$ modulo null sets. $L^{\infty+}(\mathbb{R})/_G$ is then the stalk at $G$ of this space of global sections (see \cite{PIEVIA22} for details on this example).

Cohen's forcing method applied to $\bool{MALG}$ devises a procedure which generalizes the above simultaneously to the whole topos $\bool{Sh}(\St(\bool{MALG}),\bool{CompHaus})$ (which in forcing is described as the structure $V^{\bool{MALG}}$) producing a $\bool{MALG}$-valued model of set theory. The procedure is modular and can be applied with input any complete boolean algebra $\bool{B}$ (or any extremally disconnected compact Hausdorff space), and  links the logical properties of the 
$\bool{B}$-valued structure $V^{\bool{B}}$ to the combinatorial properties of $\bool{B}$.
By suitably choosing $\bool{B}$, Cohen was able to produce a $\bool{B}$-valued model of set theory where $\bool{CH}$ gets value $0_{\bool{B}}$, or (equivalently) such that $\bool{CH}$ is false in $V^{\bool{B}}/_G$.

\section{Algebraic maximality and model companionship}\label{sec:modcomp}

It is now time to make a detour in model theory and investigate the notion of algebraic closure and its possible transfers/generalizations to other mathematical theories.\footnote{Standard model theory textbooks are \cite{CHAKEI90,HOD97,MAR02,TENZIE}.}

%
%
%
%

\subsection{Algebraic closure of structures}

Consider the signature $\bp{+,\cdot,0,1}$.
The table below  shows on the left column the algebraic theory under consideration, in the central column the first order axiomatization in the above signature, in the right column a standard model of the axioms.

\begin{center}
\begin{tabular}{| a | L | b |}
\hline
\rowcolor{LightBlue}
\textbf{Structures} & \textbf{Axioms} & \textbf{Example}\\

\hline
&&\\
 & \forall x, y\,(x\cdot y=y\cdot x) &\\
 & \forall x,y,z\, [(x\cdot y)\cdot z=x\cdot(y\cdot z)] &\\
 Commutative & \forall x \,(x\cdot 1=x\wedge 1\cdot x=x) &\\
semirings & \forall x, y\,(x + y=y + x)  &\\
with no zero&\forall x,y,z\, [(x+y)+ z=x+(y+ z)]& $\mathbb{N}$ \\
divisors&\forall y \,(x+ 0=x\wedge 0+ x=x)&\\
&\forall x,y,z\, [(x + y)\cdot z=(x\cdot y)+ (x\cdot z)]&\\ 
&\forall x, y\,[x\cdot y=0\rightarrow (x=0\vee y=0)]&\\
&&\\
\hline
&&\\
Integral && \\
domains &\forall x\exists y\, (x+y=0) & $\mathbb{Z}$ \\
&&\\
\hline

&&\\
Fields & \forall x\,[x\neq 0\rightarrow \exists y\, (x\cdot y=1)] & $\mathbb{Q}$\\
&&\\
\hline

&&\\
Algebraically  & \text{for all }n\geq 1  &\\
closed fields & \forall x_0\dots x_n\exists y\, \sum x_i\cdot y^i=0& $\mathbb{C}$\\
&&\\
\hline
\end{tabular}
\end{center}

The table describes some features of the passages from $\mathbb{N}$ to $\mathbb{Z}$, from $\mathbb{Z}$ to $\mathbb{Q}$, from $\mathbb{Q}$ to $\mathbb{C}$. 
On the semantic level these passages came along with the introduction of new operations, and the construction of structures closed under these operations (additive inverses, multiplicative inverses, closure under solutions of polynomial equations), examples of which are listed in the left-side column.
On the syntactic side the cell detected by the intersection of the first row with the central column lists the universal axioms common to all theories; the closure under the new operations of the structures on the left column is reflected by the satisfaction of corresponding $\Pi_2$-axioms, as listed in the corresponding cells of the central column.

Roughly the more closed-off a semiring with no zero-divisors is, the more $\Pi_2$-axioms it satisfies. 

Robinson came up with model theoretic notions giving an abstract description of the closed-off structures for an arbitrary universal theory.

\subsection{Existentially closed structures and model companionship}

\begin{itemize}
\item
A vocabulary $\tau$ is a list of predicates, constants and function symbols
(we are accustomed to the vocabulary
$\bp{+,\cdot,0,1,<,=}$ which we can use to write down diophantine (in)equations, such as 
$x+2y<xz+5w$). Associated to a vocabulary $\tau$ we naturally have a notion of atomic $\tau$-formula (which for 
$\bp{+,\cdot,0,1,<,=}$ are exactly the diophantine (in)equations).

\item A $\tau$-formula $\phi(x_1,\dots,x_n)$ is \emph{quantifier free} if it is a boolean combination of \emph{atomic} formulae (i.e.\ obtained by atomic formulae using conjunction $\wedge$, disjunction $\vee$, negation $\neg$, implication $\to$); for example $(x+y=1+1+z)\wedge \neg (x\cdot y<z\cdot z)$ is quantifer free for $\tau=\bp{\cdot,+,0,1,<}$.

 \item
A $\tau$-formula $\psi(x_1,\dots,x_n)$ is a $\Sigma_1$-formula if it is of the form 
\[
\exists x_0,\dots,x_k\,\phi(x_0,\dots,x_k,x_{k+1},\dots,x_n)
\] 
with $\phi(x_0,\dots,x_k,x_{k+1},\dots,x_n)$  quantifier free; for example 
\[
\exists x,y\, [(x+y=1+1+z)\wedge \neg (x\cdot y<z\cdot z)]
\] 
is $\Sigma_1$ for $\tau=\bp{\cdot,+,0,1,<}$).
 \end{itemize}

\begin{exem}
In the vocabulary $\bp{+,\cdot,0,1}$, the atomic formulae are \emph{diophantine equations} and
 the \emph{quantifier free formulae} with parameters in a ring $\mathcal{M}$ define the \emph{constructible sets} (in the sense of algebraic geometry) of $\mathcal{M}$. Below a standard example of a quantifier free  formula (with parameters in $\mathcal{M}$):
\[
\bigvee_{j=1}^{l}\qp{\bigwedge_{i=1}^{k_j} p_{ij}(a^{ij}_1,\dots,a^{ij}_{m_{ij}} ,x_1,\dots,x_n)=0\wedge  \bigwedge_{d=1}^{m_j} \neg q_{dj}(b^{dj}_1,\dots,b^{dj}_{k_{dj}},x_1,\dots,x_n)=0}
\]
with 
$a^{ij}_k,b^{dj}_k$ elements of $\mathcal{M}$ and
$p_{ij}(y_1,\dots,y_{m_{ij}} ,x_1,\dots,x_n)$, $q_{dj}(z_1,\dots,z_{k_{dj}},x_1,\dots,x_n)$ polynomials with coefficients in $\mathbb{N}$ (of degree $1$ in the $y_l,z_h$-s). 
 \end{exem}

\begin{defi}
Given a vocabulary $\tau$ and $\tau$-structures\footnote{ $\sqsubseteq$ denotes the substructure relation among the $\tau$-structures} $\mathcal{M}\sqsubseteq \mathcal{N}$,
$\mathcal{M}\prec_1\mathcal{N}$
if every $\Sigma_1$-formula
with parameters in $\mathcal{M}$ and true in $\mathcal{N}$ is true also in $\mathcal{M}$.
\end{defi}
For example
\(
\ap{\mathbb{C},+,\cdot,0,1}{\prec_1}\ap{\mathbb{C}[X],+,\cdot,0,1},
\)
but \(
\ap{\mathbb{Z},+,\cdot,0,1}{\not\prec_1}\ap{\mathbb{C},+,\cdot,0,1}\).

 \begin{defi}
 Given a $\tau$-theory $S$, a $\tau$-structure $\mathcal{M}$ is $S$-ec if:
 \begin{itemize}
 \item there is a model of $S$ $\mathcal{N}\sqsupseteq\mathcal{M}$,
 
 \item  $\mathcal{M}\prec_1 \mathcal{N}$ for any $\mathcal{N}\sqsupseteq\mathcal{M}$ which models $S$. 
 \end{itemize}
 \end{defi}
 
 \begin{exem}
 For $S$ the $\bp{+,\cdot,0,1}$-theory of \emph{integral domains} the \emph{algebraically closed fields} are exactly the $S$-ec models.
 \end{exem}

 \begin{defi}
 Given a $\tau$-theory $S$, a $\tau$-theory $T$ is the \emph{model companion} of $S$ if TFAE for any $\tau$-structure $\mathcal{M}$:
 \begin{itemize}
 \item  $\mathcal{M}$ is a model of $T$,
 
 \item  $\mathcal{M}$ is $S$-ec.
 \end{itemize}
 \end{defi}

 \begin{exem}
The $\bp{+,\cdot,0,1}$-theory of \emph{integral domains} has the $\bp{+,\cdot,0,1}$-theory of \emph{algebraically closed fields}
as its model companion.
 \end{exem}

 \begin{defi}
 A $\tau$-theory $T$ is \emph{model complete} if it is its own model companion, i.e.\ if $\mathcal{M}\prec_1\mathcal{N}$ whenever\footnote{$\mathcal{M}\prec\mathcal{N}$ holds if all formulae with parameters in $\mathcal{M}$ get the same truth value in the two structures (hence $\prec_1$ is weaker than $\prec$). Model completeness can equivalently be defined replacing $\prec_1$ by $\prec$. } $\mathcal{M}\sqsubseteq\mathcal{N}$ are models of $T$.
 \end{defi}

 \begin{exem}
The $\bp{+,\cdot,0,1}$-theory of \emph{algebraically closed fields} is \emph{model complete}.
 \end{exem}

In particular model companionship and model completeness describe a notion of algebraic closure which makes sense for an arbitrary mathematical theory, when axiomatized in a certain signature $\tau$.

Note that the signature plays a crucial role in model companionship results; for example the theory of algebraically closed fields is the model companion of the theory of integral domains in signature $\bp{+,\cdot,0,1}$.\footnote{It is also the model companion of the theory of fields in this same signature.} It is not anymore so in the signature $\bp{+,\cdot,^{-1},0,1}$ where to interpret $^{-1}$ in an integral domain we use the axiom:
\[
\forall x\,[\exists y \,(x\cdot y=1\wedge x\cdot x^{-1}=1)\vee(\neg\exists y\, (x\cdot y=1)\wedge x^{-1}=0)].
\]
With this axiom we still get that
\[
\mathcal{M}=\ap{\mathbb{C},+,\cdot,^{-1},0,1}\mathrel{\sqsubseteq}\mathcal{N}=\ap{\mathbb{C}[X],+,\cdot,^{-1},0,1},
\]
but not that
\(
\mathcal{M}\mathrel{\prec_1}\mathcal{N}:
\)
$\exists x\,(x\neq 0\wedge x^{-1}=0)$ is true in $\mathcal{N}$ (as witnessed by $X$), but false in $\mathcal{M}$.

In particular $\ap{\mathbb{C},+,\cdot,^{-1},0,1}$ is not anymore $T$-ec for $T$ the theory of integral domains in signature $\bp{+,\cdot,^{-1},0,1}$. 
Note however that the theory of algebraically closed fields is still model complete\footnote{It is also still the model companion of the theory of fields in $\bp{+,\cdot,^{-1},0,1}$.} even when formalized in signature $\bp{+,\cdot,^{-1},0,1}$.

\section{Algebraic closure for set theory} \label{sec:algclsetth}
We now want to use the notion of model companionship to describe what is ``algebraic closure'' for set theory. This will allow to formulate a reasonably accessible definition of Woodin's axiom $(*)$.
\subsection{The right vocabulary for set theory}
The standard axiomatization of set theory in textbooks is done in vocabulary $\bp{\in}$, eventually with extra symbol $\subseteq$. However this is not giving an efficient formalization.
For example in the $\bp{\in}$-vocabulary the notion of ordered pair is formalized by means of
\emph{Kuratowski's trick:} $\ap{y,z}$ is (coded by) the set $\bp{\bp{y},\bp{y,z}}$. 

The standard $\in$-formula expressing $x=\ap{y,z}$ is
\[
\exists t\exists u\;[\forall w\,(w\in x\leftrightarrow w=t\vee w=u)
\wedge\forall v\,(v\in t\leftrightarrow v=y)
\wedge\forall v\,(v\in u\leftrightarrow v=y\vee v=z)].
\]
We do not regard the concept of ordered pair as a complicated one. It should then be formalizable by a simple formula. Accordingly many other basic set theoretic concepts such as that of being an $n$-ary relation, a function, the emptyset, the set of natural numbers,\dots should all be formalizable by simple formulae.

The lightface $\Delta_0$-properties isolate the simplest set theoretic properties and include the following classes:
\begin{itemize}
\item $\bp{R\in V: R\text{ is an $n$-ary relation}}$,
\item $\bp{f\in V: f\text{ is a function}}$,
\item $\bp{\ap{a,b}: a\subseteq b\text{ are sets}}$,
\item\dots
\end{itemize}
\begin{defi}\parencite[Def. IV.3.5]{KUNEN}

An $\bp{\in}$-formula $\phi$ is a $\Delta_0$-formula if all its quantified variables are bounded to range in a set (e.g.\ $y\subseteq z\, \equiv \,\forall x\,(x\in y\rightarrow x\in z)\, \equiv\, \forall x\in y\,(x\in z)$ is a $\Delta_0$-formula).

A lightface $\Delta_0$-property is a class of sets whose extension is given by a $\Delta_0$-formula.
\end{defi}


It is thus natural to expand the vocabulary of set theory as follows:

\begin{defi}
\emph{}

$\in_{\Delta_0}$ has  the following list of symbols:
\begin{itemize} 
\item
constants for $\emptyset,\mathbb{N}$, 
\item
relation symbols $R_\phi$ for any $\Delta_0$-formula $\phi(x_1,\dots,x_n)$, 
\item
function symbols $G_i$ for the list of ten {G}\"{o}del operations given in \cite[Def. 13.6]{JECHST}.
\end{itemize}
\end{defi}
The base $\in_{\Delta_0}$-theory for sets interprets the new symbols according to their expected meaning, for example it must include the following axioms:
\begin{itemize}
\item $\forall \vec{x}\,[R_\phi(\vec{x})\leftrightarrow\phi(\vec{x})]$ whenever $\phi$ is a $\Delta_0$-formula
\item $\forall\vec{x},y\, (R_{\phi_i}(\vec{x},y)\leftrightarrow G_i(\vec{x})=y)$ if $\phi_i(\vec{x},y)$ is the $\Delta_0$-formula describing the graph of the  {G}\"{o}del operation $G_i$.
\item $\forall x\,(x\not\in \emptyset)$.
\item $\forall x\, (x\in\mathbb{N}\leftrightarrow x\text{ is a finite Von Neumann ordinal})$. 
\end{itemize}
Note that
$(x\text{ is a finite Von Neumann ordinal})$ is (formalizable by) a $\Delta_0$-formula.

We now give an explicit first order axiomatization of $\bool{MK}$ in signature $\in_{\Delta_0}\cup\bp{\bool{Set},V}$ where $\bool{Set}$ is a unary predicate symbol for the property of being a set and $V$ is a constant which denotes the universe of all sets.
We stipulate for the sake of readability that smallcase letters indicate sets, uppercase letters indicate classes.
To the above list of $\in_{\Delta_0}$-axioms\footnote{To be picky, the above list of axioms should be reformulated so that quantifiers apply only to variables ranging over sets.} we add:

%
%
%
%
%
%
%
%
%
%
%

\begin{description}
\item[Universal axioms]
\emph{}
\vspace{0.1cm}
\begin{description}
\item[Extensionality] $\forall X,Y \,[(X\subseteq Y\wedge Y\subseteq X)\leftrightarrow X=Y]$.

\item[Comprehension (a)] $\forall X\, (\bool{Set}(X)\leftrightarrow X \in V)\wedge\forall X\,(X\subseteq V)$.

\item[Foundation] 
\[
\forall F\,[(F\text{ is a function }\wedge\dom(F)=\mathbb{N})\,\rightarrow \exists n\in \mathbb{N} \,F(n+1)\not\in F(n)].
\]

\end{description}

\vspace{0.2cm}

\item[Existence Axioms]
\emph{}
\vspace{0.1cm}
\begin{description}

\item[Emptyset] $\bool{Set}(\emptyset)$, 
\item[Infinity] $\bool{Set}(\mathbb{N})$. 
\end{description}

\vspace{0.2cm}

\item[Basic construction principles]
\emph{}
\vspace{0.1cm}
\begin{description}
\item[{G}\"{o}del operations] For each $i=1,\dots,10$ 
\[
\forall x_0,\dots,x_{k_i}[(\bigwedge_{i=0}^k(x_i\in V)\rightarrow \exists! z\,(z\in V\wedge R_{\phi_i}(\vec{x},z))],
\]
 where $\phi_i$ is the $\Delta_0$-formula whose extension is the graph of the {G}\"{o}del operation $G_i$.
\item[Separation] $\forall P,x\, [\bool{Set}(x)\rightarrow \bool{Set}(P\cap x)]$.
\end{description}

\vspace{0.2cm}

\item[Strong construction principles]
\emph{}
\vspace{0.1cm}
\begin{description}
\item[Comprehension (b)] For every $\in_{\Delta_0}$-formula $\psi(x_0,\dots,x_n,\vec{Y})$ 
\[
\forall \vec{Y}\,\exists Z\,\forall x\,[x\in Z\leftrightarrow (x\in V\wedge \exists x_0,\dots,x_n (x=\ap{x_0,\dots,x_n}\wedge \psi(x_0,\dots,x_n,\vec{Y})))].
\]
\item[Replacement]
\[
\forall F, x \, [(F \text{ is a function } \wedge \bool{Set}(x)\wedge (x\subseteq\dom(F)))\rightarrow \bool{Set}(F[x])].
\]

\item[Powerset] 
\[
\forall x\,[\bool{Set}(x)\rightarrow \exists ! y\, [\bool{Set}(y)\wedge (\forall z\,(z\in y\leftrightarrow z\subseteq x))]].
\]

\item[Global Choice]

{\small
\begin{align*}
\forall F[&&\\
&F\text{ is a function}\wedge \forall x\,(x\in\dom(F)\rightarrow F(x)\neq\emptyset)&\\
\rightarrow&&\\ 
&\exists G\, (G\text{ is a function}\wedge \dom(G)=\dom(F)\wedge 
\forall x\,(x\in\dom(G)\rightarrow G(x)\in F(x))&\\
]&&
\end{align*}
}
\end{description}
\end{description}

An $\in_{\Delta_0}$ model of $\MK$ is a ``sorted'' structure $(\mathcal{C},V,\in_{\Delta_0})$, where 
$\mathcal{C}$ is the family of all classes, and $V$ is an element of $\mathcal{C}$ whose extension is the subfamily of $\mathcal{C}$ given by sets. 
One can also check that in this case
$(V,\in_{\Delta_0})$ is a model of the $\bool{ZFC}$-axioms (the classical presentation of set theory which avoids the mention of proper classes).

With a certain degree of approximation, it is customary to denote a model  $(\mathcal{C},V,\in_{\Delta_0})$
by its second component $V$, as $\mathcal{C}$ could be recovered as the family of subclasses of $V$.
However this is not entirely correct as there can be $(\mathcal{C},V,\in_{\Delta_0})$, $(\mathcal{D},V,\in_{\Delta_0})$ both models of $\bool{MK}$ with $\mathcal{C}\neq\mathcal{D}$. When confusion on this issue may arise we will be more explicit. In general we stick to the convention to denote a model of (some of the axioms of)
$\bool{MK}$ by $(\mathcal{C},V,\in_{\Delta_0})$ and when we write just $(V,\in_{\Delta_0})$ we are considering the substructure whose elements are sets and not proper classes.

\subsection{The $H_\kappa$s}

A finite set may not be simple, for example to understand the singleton $\bp{\mathbb{R}}$ we need to know $\mathbb{R}$. We want to stratify sets according to the cardinalities required to generate them.
Recall (Def. \ref{def:Hkappa} and Def. \ref{def:Haleph0}) that a set $X$ is \emph{hereditarily of size less than (at most) $\kappa$} if $\trcl(X)$ has size less than (at most) $\kappa$.
%
%
%
%

\begin{itemize}
\item $\bp{\mathbb{R}}$ is not hereditarily countable (i.e.\ $\bp{\mathbb{R}}\notin H_{\aleph_1}$);
\item Any subset $A$ of $\mathbb{N}$ is hereditarily countable (i.e.\ $A\in H_{\aleph_1}$);
\item
$\mathbb{Q}$ and $\mathbb{Z}$ as defined in any textbook are hereditarily countable (i.e.\ in $H_{\aleph_1}$);
\item
$\mathbb{R}$ and $\pow{\mathbb{N}}$ are subsets of $H_{\aleph_1}$ (but not elements!);
\item $\pow{\mathbb{N}}$ is definable by the atomic $\in_{\Delta_0}$-formula 
$(x\subseteq\mathbb{N})$ in the structure $\ap{H_{\aleph_1},\in_{\Delta_0}}$;
\item similarly for $\mathbb{R}$ or for (a representative of the homeomorphism class of) any Polish space.
\end{itemize}

\begin{itemize}
\item $\pow{\aleph_1}$ is definable by the atomic $\in_{\Delta_0}$-formula $(x\subseteq\aleph_1)$ in parameter $\aleph_1$ (the first uncountable ordinal)
 in the structure $\ap{H_{\aleph_2},\in_{\Delta_0}}$,
\item $\NS$, the non-stationary ideal on $\aleph_1$, is $\Sigma_1$-definable in parameter $\aleph_1$ in the same structure (see Section \ref{sec:wooax*} for a definition of $\NS$).
\end{itemize}

\[
H_{\aleph_0}\subseteq H_{\aleph_1}\subseteq H_{\aleph_2}\subseteq\dots\subseteq H_{\kappa^+}\subseteq\dots
\]

\[
V=\bigcup\bp{H_{\lambda}: \lambda\text{ an infinite cardinal}}
\]

$(\pow{H_{\aleph_i}},H_{\aleph_i},\in_{\Delta_0})$ for $i=1,2$ are models of all axioms of $\MK$ with the exception of powerset. 
The role of $V$ in these models is played by $H_{\aleph_i}$.\footnote{Similarly it can be shown that for $i=1,2$, $(H_{\aleph_i},\in_{\Delta_0})$ are model of all axioms of $\ZFC$ with the exception of Powerset.}

$(\pow{H_\kappa},H_{\kappa},\in_{\Delta_0})$ models $\MK$ if and only if $\kappa$ is inaccessible.

\subsection{Existentially closed structures for set theory}
\begin{theo}[Levy absoluteness]
Let $\kappa$ be an infinite cardinal.

Then
\[
\ap{H_{\kappa^+},\in_{\Delta_0},A: A\subseteq \pow{\kappa}}\prec_1\ap{V,\in_{\Delta_0},A: A\subseteq\pow{\kappa}}
\]
\end{theo}

We argue below that $H_{\aleph_1}$ and $H_{\aleph_2}$ provide natural candidates for existentially closed structures for set theory. The choice of whether to focus on $H_{\aleph_1}$ or $H_{\aleph_2}$ depends on the signature in which one formalizes set theory.

\section{Algebraic maximality for $\pow{\mathbb{N}}$ and generic absoluteness} \label{sec:algclsont}

We connect here the notion of algebraic closure given by existentially closed model, to generic absoluteness results for second order arithmetic.

\begin{theo}[Shoenfield, 1961]
Let $V[G]$ be a forcing extension of $V$. Then\footnote{Note that $\ap{H_{\aleph_1}^{V[G]},\in_{\Delta_0}}\prec_1\ap{V[G],\in_{\Delta_0}}$ follows from Levy Absoluteness applied in $V[G]$.}
\[
\ap{H_{\aleph_1},\in_{\Delta_0}}\prec_1
\ap{H_{\aleph_1}^{V[G]},\in_{\Delta_0}}.
\]
\end{theo}

A  key observation is that $\pow{\mathbb{N}}$ can be identified with Cantor's space $2^\mathbb{N}$ via characteristic functions. With this identification the universally Baire subsets\footnote{Def. \ref{def:UB} was given only on $\mathbb{R}^k$ but makes sense for any uncountable Polish space $X$ (a second countable topological space whose topology can be induced by a complete metric) which is locally compact. Actually the families of universally Baire subsets of uncountable locally compact Polish spaces $X$ and $Y$ can be identified modulo any Borel isomorphism existing between the two spaces \parencite{FENMAGWOO,kechris:descriptive}.} of $2^{\mathbb{N}}$ describe a particularly nice $\sigma$-algebra contained in $\pow{\pow{\mathbb{N}}}$. 

$\bool{UB}^V$ denotes the family of universally Baire subsets of $2^\mathbb{N}$ existing in $V$.
Given $V[G]$ forcing extension of $V$,
each universally Baire set $A$ in $V$ has a canonical extension $A^{V[G]}$ to a universally Baire set of 
$V[G]$ \parencite{FENMAGWOO}.
%

The following is an elaboration which rephrases in model theoretic terminology deep results of Feng, Magidor, Woodin, Steel, Martin (see \cite{VIAVENMODCOMP}):
\begin{theo}[Woodin 1985 + Steel, Martin, 1989 + Feng, Magidor, Woodin, 1992 +
  V., Venturi, 2019]
Assume there is a proper class of Woodin's cardinals.
Then the theory of the structure 
\[
\ap{H_{\aleph_1}^{V},\in_{\Delta_0},A^{V}: { A\in \bool{UB}^V}}
\]
is model complete and is the model companion of the theory of 
\[
\ap{V[H],\in_{\Delta_0},{ A^{V[H]}: A\in\bool{UB}^V}}
\] 
whenever $V[H]$ is some generic extension of $V$.
\end{theo}

\subsection{Algebraic maximality for $\pow{\mathbb{N}}$: a summary}
The table below summarizes the effects of large cardinals on the algebraic closure properties of set theory relative to $\pow{\mathbb{N}}$ (or -better- to $H_{\aleph_1}$):

\vspace{0.3 cm}

\begin{tabular}{| a | L |}
\hline
\rowcolor{LightBlue}
\textbf{Theory} & \textbf{Degree of algebraic closure} \\
\hline
& \ap{H_{\aleph_1},\in_{\Delta_0},{ A:  A\in \bool{UB}^V}}\\
$\bool{MK}$ &\text{ is }\Sigma_1\text{-elementary in}\\
&\ap{V{[G]},\in_{\Delta_0},{ A^{V[G]}: A\in\bool{UB}^V}}\\
&\text{ for all generic extension }V{[G]}\text{ of }V\\
\hline
&\\
$\bool{MK+}$ &
\ap{H_{\aleph_1},\in_{\Delta_0}, { A: \, A\in \bool{UB}^V}}\\
large cardinals&\text{has a \textbf{model complete theory}}\\
 &\text{which is the \textbf{model companion} of the theory of}\\
&\ap{V{[G]},\in_{\Delta_0},{ A^{V[G]}: A\in\bool{UB}^V}}\\
&\text{ for all generic extension }V{[G]}\text{ of }V\\
\hline
\end{tabular}
%
%
%
%
%
%

\section{Algebraic maximality for $\pow{\aleph_1}$ and axiom $(*)$} \label{sec:wooax*}

The notions of club and stationarity are given for the sake of completeness. One could take the extension of these concepts as a blackbox and still get some useful insights on Woodin's axiom $(*)$.
\begin{itemize}
\item
$C$ is a club subset\footnote{We here represent $\aleph_1$ by the least Von Neumann ordinal in the equivalence class of $\aleph_1$.} of $\aleph_1$ if $\cup C=\aleph_1$ and 
for all $\beta\notin C$ there is $\alpha\in\beta$ such that $(\alpha,\beta]\cap C$ is empty (where 
 $(\alpha,\beta]$ is given by those ordinals $\gamma$ with $\alpha\in\gamma$ and $\gamma\in\beta$ or $\gamma=\beta$).

\item
$S\subseteq\aleph_1$ is stationary if for all $C$ club subset of $\aleph_1$ $S\cap C$ is non-empty.

\item $\NS\subseteq\pow{\aleph_1}$ is the ideal of non-stationary subsets of $\aleph_1$ (i.e.\ subsets
disjoint from some club).

\end{itemize}
%
%

\begin{defi}\label{def:SSP}
Let $\bool{B}$ be a complete boolean algebra. 
$\bool{B}$ is $\bool{SSP}$ if whenever $V[G]$ is a forcing extension of $V$ by $\bool{B}$
\[
\ap{H_{\aleph_2},\in_{\Delta_0},\NS^V}{\sqsubseteq}
\ap{V[G],\in_{\Delta_0},\NS^{V[G]}}.
\]
\end{defi}

\begin{defi}[\cite{BAG00}]
\textbf{Bounded Martin's maximum} $\bool{BMM}$ holds if whenever $\bool{B}$ is an $\bool{SSP}$ complete boolean algebra and 
$V[G]$ is a forcing extension of $V$ by $\bool{B}$
\[
\ap{H_{\aleph_2},\in_{\Delta_0}}{\prec_1}
\ap{H_{\aleph_2}^{V[G]},\in_{\Delta_0}}.
\]
\end{defi}
Compare it with Shoenfield's theorem and Levy absoluteness.

\begin{theo}[\cite{BAG00}]
$\MM$ implies $\bool{BMM}$.
\end{theo}

\begin{defi}[\cite{woodinBOOK}]
 $\bool{BMM}^{++}$ holds if whenever $\bool{B}$ is an $\SSP$ complete boolean algebra and 
$V[G]$ is a forcing extension of $V$ by $\bool{B}$
\[
\ap{H_{\aleph_2},\in_{\Delta_0},\NS^V}{\prec_1}
\ap{H_{\aleph_2}^{V[G]},\in_{\Delta_0},\NS^{V[G]}}.
\]
\end{defi}

\begin{theo}[\cite{woodinBOOK}]
$\MM^{++}$ implies $\bool{BMM}^{++}$.
\end{theo}

\subsection{Applications of $\bool{BMM}$}

The consequences of bounded forcing axioms are almost the same of those of $\bool{MM}$ although in some cases new proofs had to be found.

Assume $\bool{BMM}$. Then:
\begin{itemize}
\item
 $2^{\aleph_0}=\aleph_2=\aleph_1^+$ \parencite{TOD02}. 

\item
Whitehead's conjecture on free groups is false.

\item Kaplansky's conjecture on Banach algebras holds.


\end{itemize}
On the other hand Moore's result on the existence of a five element basis for uncountable linear orders and Farah's result establishing that $\MM$ implies all automorphism of the Calkin algebra are inner are not known to follow from $\bool{BMM}^{++}$.

\subsection{Woodin's axiom $(*)$}

We let 
$\bool{UB}^V$ denote the family of universally Baire subsets of $\pow{\mathbb{N}}$ existing in $V$.

\begin{rema} 
Let $\bool{B}$ be an $\SSP$ complete boolean algebra and 
$V[G]$ be a forcing extension of $V$ by $\bool{B}$. Then
{\small
\[
\ap{H_{\aleph_2},\in_{\Delta_0},\NS, { A:\, A\in \bool{UB}^V}}{\sqsubseteq}
\ap{H_{\aleph_2}^{V[G]},\in_{\Delta_0},\NS^{V[G]}, { A^{V[G]}: A\in \bool{UB}^V}}.
\]
}
\end{rema}

It is not clear to me whether the following definition is due to Woodin or Shelah or Goldstern:
\begin{defi}[Woodin, Shelah, Goldstern?]
$\bool{UB}\text{-}\bool{BMM}^{++}$ holds if whenever $\bool{B}$ is an $\bool{SSP}$ complete boolean algebra and 
$V[G]$ is a forcing extension of $V$ by $\bool{B}$
{\small
\[
\ap{H_{\aleph_2},\in_{\Delta_0},\NS, { A:\, A\in \bool{UB}^V}}{\prec_1}
\ap{H_{\aleph_2}^{V[G]},\in_{\Delta_0},\NS^{V[G]}, { A^{V[G]}: A\in \bool{UB}^V}}.
\]
}
\end{defi}

By variations of the results in \textcite{woodinBOOK}, one gets:

\begin{theo}[Woodin]
$\MM^{++}$ implies $\bool{UB}\text{-}\bool{BMM}^{++}$.
\end{theo}

$(*)_{\bool{UB}}$ is a natural strengthening of Woodin's axiom $(*)$.
\begin{theo}[\cite{ASPSCH(*)}]
Assume there is a proper class of Woodin cardinals.
Then
$(*)_{\bool{UB}}$ if and only if $\bool{UB}\text{-}\bool{BMM}^{++}$.
\end{theo}

If one is interested in Woodin's axiom $(*)$, here is an equivalent reformulation of it:
\begin{theo}[\cite{ASPSCH(*)}]
Assume there is a proper class of Woodin cardinals.
Then
Woodin's axiom $(*)$ holds if and only if 
whenever $\bool{B}$ is an $\SSP$ complete boolean algebra and 
$V[G]$ is a forcing extension of $V$ by $\bool{B}$
\[
\ap{H_{\aleph_2},\in_{\Delta_0},\NS, { A:\, A\text{ is in }\pow{\mathbb{R}}^{L(\mathbb{R})^V}}}
\]
is $\Sigma_1$-elementary in
\[
\ap{H_{\aleph_2}^{V[G]},\in_{\Delta_0},\NS^{V[G]},{ A^{V[G]}:  A\text{ is in }\pow{\mathbb{R}}^{L(\mathbb{R})^V}}}.
\]
\end{theo}

The original formulation of $(*)$ can be found in \textcite{woodinBOOK}. $(*)$ follows from $(*)_{\bool{UB}}$
once one notes that any set of reals definable in $L(\mathbb{R})$ is universally Baire (assuming the existence of a proper class of Woodin cardinals).

\subsection{Woodin's axiom $(*)$ and model completeness for the theory of $H_{\aleph_2}$}

Recall that $\psi$ is a $\Pi_2$-sentence if it is of the form ${\forall\vec{x}}\,{\exists\vec{y}}\phi({\vec{x}},{\vec{y}})$ with $\phi(\vec{x},\vec{y})$ quantifier free.

In signature $\in_{\Delta_0}$ $\neg\bool{CH}$ can be formalized by the $\Pi_2$-sentence in parameter 
$\aleph_1$ (the first uncountable ordinal/cardinal):
\[
\forall f\,[(\underbrace{f\text{ is a function}}_{\Delta_0(f)}\wedge 
\underbrace{\dom({f})={\aleph_1}}_{\Delta_0(f,\aleph_1)})\rightarrow 
{\exists r}\,(\underbrace{{r}\subseteq\mathbb{N}}_{\Delta_0(r,\mathbb{N})}\wedge \underbrace{{r}\not\in \ran({f})}_{\Delta_0(r,f)})]
\]

Recall that $\NS$ is saturated if the boolean algebra $\pow{\aleph_1}/_{\NS}$ has only partitions of size at most $\aleph_1$.

We quote the following facts on the non-stationary ideal:

\begin{itemize}
\item
Assume $\NS$ is saturated. Then it is precipitous.
\item
Assume $\MM$. Then $\NS$ is saturated \parencite{FORMAGSHE}.
\item
$\NS$ is precipitous is consistent with $\bool{CH}$.
\end{itemize}

\begin{theo}[\cite{woodinBOOK}]
Assume there is a proper class of supercompact cardinals, \emph{ Sealing\footnote{See Def. \ref{def:sealing}.}, and $\NS$ is precipitous}.
TFAE:
\begin{itemize}
\item $(*)_{\bool{UB}}$ (or $\bool{UB}\text{-}\bool{BMM}^{++}$).

\item For any  $\Pi_2$-sentence\footnote{Among which $\neg\bool{CH}$ and a strong form of $2^{\aleph_0}=\aleph_2$.} $\psi$ for $\in_{\Delta_0}\cup\bp{\aleph_1,\NS}\cup\bp{A:A\in \bool{UB}^V}$
\[
\ap{H_{\aleph_2},\in_{\Delta_0},\aleph_1,\NS, A:\, A\in {\bool{UB}^V}}\models\psi
\]
\centerline{if and only if}

\centerline{
\textbf{$\psi$ is true in $H_{\aleph_2}^{V[G]}$ for some forcing extension $V[G]$ of $V$}.}
\end{itemize}
\end{theo}
 Sealing can be removed if one replaces $\bool{UB}^V$ with $\pow{\mathbb{R}}^{L(\Ord^\mathbb{N})}$ in the formulation of $\bool{UB}\text{-}\bool{BMM}^{++}$  and in all relevant spots (note that the sets of reals definable in 
 $L(\Ord^\mathbb{N})$ are universally Baire assuming the existence of class many supercompact cardinals).
The large cardinal assumptions are far stronger than needed.

\begin{theo}[\cite{viale2021absolute}]
Assume there is a proper class of supercompact cardinals, \emph{{ Sealing}, and $\NS$ is precipitous}.
TFAE:

\begin{itemize}
\item $(*)_{\bool{UB}}$ (or $\bool{UB}\text{-}\bool{BMM}^{++}$).

\item The theory $T$ of the structure 
\[
\mathcal{M}=\ap{H_{\aleph_2},\in_{\Delta_0},\aleph_1,\NS, A:\, A\in {\bool{UB}^V}}
\]
is the \textbf{model companion} of the theory $S$ of the structure 
\[
\ap{V,\in_{\Delta_0},\aleph_1,\NS, A:\, A\in {\bool{UB}^V}}.
\]
(i.e.\ $T$ is \textbf{model complete})

\item Letting $S_{\forall\vee\exists}$ be the boolean combination of existential sentences
which are in $S$, and $\psi$ be a $\Pi_2$-sentence,

\centerline{\textbf{$\mathcal{M}$ models $\psi$ if and only $\psi+S_{\forall\vee\exists}$ is consistent}.}

\item For any  $\Pi_2$-sentences $\psi$ 
\[
\ap{H_{\aleph_2},\in_{\Delta_0},\NS, A:\, A\in {\bool{UB}^V}}\models\psi
\]
\centerline{if and only if}

\centerline{
\textbf{$\psi$ is true in $H_{\aleph_2}^{V[G]}$ for some forcing extension $V[G]$ of $V$}.}

\end{itemize}
\end{theo}

 { Sealing} can be removed if one replaces $\bool{UB}^V$ with $\pow{R}^{L(\Ord^\mathbb{N})}$ in the formulation of $\bool{UB}\text{-}\bool{BMM}^{++}$ and in the relevant spots.
 
 \subsection{Sealing}
 For the sake of completeness a form of Sealing sufficient to prove both theorems is the following:

 Given $(\mathcal{D},W,\in_{\Delta_0})$ transitive model of $\bool{MK}$, let 
 $N^W$ be the set $\pow{H_{\aleph_1}}^{L(\bool{UB})^W}$, where $L(\bool{UB})^W$ is the smallest transitive model of $\ZF$ containing $\bool{UB}^W$.

 \begin{defi}[Woodin]\label{def:sealing}
\textbf{Sealing} holds in a model $(\mathcal{C},V,\in_{\Delta_0})$ of $\bool{MK}$ if:
\begin{itemize}
\item
there are class many Woodin cardinals;
\item
the theory $T$
of $(N^V,H_{\aleph_1}^{V},\in_{\Delta_0})$ is model complete;
\item
$(N^{V[G]},H_{\aleph_1}^{V[G]},\in_{\Delta_0})$ models $T$ whenever $V[G]$ is a generic extension of $V$.
\end{itemize}
%
 \end{defi}
 
 It is not known whether Sealing follows rightaway from large cardinals; however
 Woodin has extablished a strong form of consistency for sealing relative to large cardinals:
 \begin{theo}[Woodin] 
 Assume $V$ models $\kappa$ is supercompact and there is a proper class of Woodin cardinals. Let $V[H]$ be a generic extension of $V$ where 
 $\kappa$ is countable. Then sealing holds in $V[H]$.
 \end{theo}

In particular assuming large cardinals, Sealing can be forced to be true, and once it is true in a forcing extension there is no way to force it to become false in a further forcing extension. See \cite[Section 3.4]{STATLARSON}
for details.

\subsection{Algebraic maximality for $\pow{\aleph_1}$: a summary}
The table below summarizes the effects of large cardinals and forcing axioms on the algebraic closure properties of set theory relative to $\pow{\aleph_1}$ (or -better- to $H_{\aleph_2}$):

\vspace{0.3 cm}

\begin{tabular}{| a | L |}
\hline
\rowcolor{LightBlue}
\textbf{Theory} & \textbf{Degree of algebraic closure} \\
\hline
&\\
& \ap{H_{\aleph_2}^V,\in_{\Delta_0},\aleph_1^V,{ \NS},{ A^V:  A\in \bool{UB}^V}}\\
$\bool{MK}$ &\text{is  a \emph{substructure} of }\\
&\ap{V{[G]},\in_{\Delta_0},\aleph_1^{V[G]},{ \NS^{V[G]}},{ A^{V[G]}: A\in\bool{UB}^V}}\\
&\text{ for all generic extension }V{[G]}\text{ of }V\text{ by an }{\bool{SSP}}\text{-forcing}\\
&\\
\hline
&\\
& \ap{H_{\aleph_2}^V,\in_{\Delta_0},\aleph_1^V,{ \NS^V},{ A^V:  A\in \bool{UB}^V}}\\
$\bool{MK}+$ &\text{is  a }\Sigma_1\text{\emph{-substructure} of}\\
forcing&\ap{V{[G]},\in_{\Delta_0},\aleph_1^{V[G]},{ \NS^{V[G]}},{ A^{V[G]}: A\in\bool{UB}^V}}\\
axioms&\text{ for all generic extension }V{[G]}\text{ of }V\text{ by an }{ \bool{SSP}}\text{-forcing}\\
&\\
\hline
&\\
$\bool{MK+}$& \text{for all generic extension }V{[G]}\text{ of }V\text{ the theories of}\\
large cardinal &\ap{V{[G]},\in_{\Delta_0},\aleph_1^{V[G]},{ \NS^{V[G]}}, { A^{V[G]}: \, A\in \bool{UB}^V}}\\
axioms &\text{have the same \textbf{model companion} theory}\\
&\\
\hline
&\\
& \text{for all generic extension }V{[G]}\text{ of }V\text{ the theories of}\\
$\bool{MK+}$&
\ap{V{[G]},\in_{\Delta_0},\aleph_1^{V[G]},{\NS^{V[G]}}, { A^{V[G]}: \, A\in \bool{UB}^V}}\\
large cardinals $+$&\text{have as \textbf{model companion} the theory of}\\
forcing &\ap{H_{\aleph_2}^V,\in_{\Delta_0},\aleph_1^V,{ \NS^V},{ A^V:  A\in \bool{UB}^V}}\\
axioms&\\
&\\
\hline
\end{tabular}

\printbibliography

\end{document}


%% file: macros.tex
\newcommand{\Ord}{\ensuremath{\mathrm{Ord}}}

\newcommand{\ZFC}{\ensuremath{\mathsf{ZFC}}}
\newcommand{\ZF}{\ensuremath{\mathsf{ZF}}}
\newcommand{\MK}{\ensuremath{\mathsf{MK}}}

\DeclareMathOperator{\dom}{dom}
\DeclareMathOperator{\ran}{ran}

\DeclareMathOperator{\trcl}{trcl}

\DeclareMathOperator{\St}{St}

\newcommand{\NS}{\ensuremath{\mathsf{NS}}} 

\newcommand{\bool}[1]{\mathsf{#1}}

\newcommand{\pow}[1]{\mathcal{P}\left(#1\right)}

\newcommand{\qp}[1]{\left[ #1 \right]}
\newcommand{\Qp}[1]{\left\llbracket #1 \right\rrbracket}
\newcommand{\ap}[1]{\langle #1 \rangle}
\newcommand{\bp}[1]{\left\lbrace #1 \right\rbrace}

\newcommand{\MM}{\ensuremath{\text{{\sf MM}}}} 

\newcommand{\SSP}{\ensuremath{\text{{\sf SSP}}}}

\newcommand{\FA}{\ensuremath{\text{{\sf FA}}}}
